\documentclass[11pt]{amsart}

\newtheorem{theorem}{Theorem}[section]
\newtheorem{proposition}{Proposition}[section]
\newtheorem{lemma}{Lemma}[section]
\newtheorem{corollary}{Corollary}[section]

\theoremstyle{definition}

\numberwithin{equation}{section}

\newcommand{\Nat}{{\mathbb N}}

\newcommand{\Real}{{\mathbb R}}

\newcommand{\B}{{\mathcal B}}
\newcommand{\Bnat}{{\mathcal B^{nat}}}
\newcommand{\hilbert}{{\mathcal H}}
\newcommand{\kspace}{{\mathcal K}}
\newcommand{\M}{{\textbf{M}}}

\def\hlimit{\overset{\hilbert}{\rightarrow}}

\begin{document}

\title[On Nyman-Beurling Criterion]{A strengthening of the Nyman-Beurling criterion for the Riemann Hypothesis, 2}

\author{Luis B\'{a}ez-Duarte}

\date{28 March 2002}

\email{lbaez@ccs.internet.ve}


\begin{abstract}
Let $\rho(x)=x-[x]$, $\chi=\chi_{(0,1)}$. In $L_2(0,\infty)$ consider the subspace $\B$ generated by $\{\rho_a|a\geq1\}$
where $\rho_a(x):=\rho\left(\frac{1}{ax}\right)$. By the Nyman-Beurling criterion the Riemann hypothesis is equivalent
to the statement $\chi\in\overline{\B}$. For some time it has been conjectured, and proved in the first version of this
paper, posted in arXiv:math.NT/0202141 v2, that the Riemann hypothesis is equivalent to the stronger statement that 
$\chi\in\overline{\Bnat}$ where $\Bnat$ is the much smaller subspace generated by $\{\rho_a|a\in\Nat\}$. This second
version differs from the first in showing that under the Riemann hypothesis the distance between $\chi$ and
$-\sum_{a=1}^n\mu(a)e^{-c\frac{\log a}{\log\log n}}\rho_a$ is of order $(\log\log n)^{-1/3}$.
\end{abstract}

\maketitle

\section{Introduction}
We denote the fractional part of $x$ by  $\rho(x)=x-[x]$, and let $\chi$ stand for the characteristic function of the
interval $(0,1]$. $\mu$ denotes the M\"{o}bius function. We shall be working in the Hilbert space

$$
\hilbert:=L_2(0,\infty),
$$
where our main object of interest is the subspace of \textit{Beurling functions}, defined
as the linear hull of the family $\{\rho_a|1\leq a\in\Real\}$ with

$$
\rho_a(x):=\rho\left(\frac{1}{ax}\right).
$$
The much smaller subspace $\Bnat$ of \textit{natural Beurling functions} is generated by $\{\rho_a|a\in\Nat\}$. The
Nyman-Beurling criterion (\cite{nyman}, \cite{beurling}) states, in a slightly modified form
\cite{notes3} (the original formulation is related to $L_2(0,1)$), that the Riemann hypothesis is equivalent to the
statement that 

$$
\chi\in\overline{\B},
$$
\ \\
but it has recently been conjectured by several authors\footnote{see \cite{BR1}, \cite{IUO}, \cite{AB}, \cite{notes3},
\cite{notes1}, \cite{conrey}, \cite{frankenhuysen},
\cite{landreau}, \cite{lee}, \cite{nikolski}, \cite{vasyunin}, \cite{vasyunin2}} that this condition
could be substituted by $\chi\in\overline{\Bnat}.$ We state this as a theorem to be proved below.

\begin{theorem}\label{natconj}
The Riemann hypothesis is equivalent to the statement that

$$
\chi\in\overline{\Bnat}.
$$
\end{theorem}
\ \\
To properly gauge the strength of this theorem note this: not only is $\Bnat$ a rather
thin subspace of $\B$ but, as is easily seen, it is also true that $\overline{\B}$ is much larger than
$\overline{\Bnat}$. 

By necessity all authors have been led in one way or another to the \textit{natural approximation}
\begin{equation}\label{natapprox}
F_n:=\sum_{a=1}^n\mu(a)\rho_a,
\end{equation}
which tends to $-\chi$ both a.e. and in $L_1$ norm when restricted to $(0,1)$ (see
\cite{BR1}), but which has been shown (\cite{IUO}, \cite{AB}) to diverge in $\hilbert$. In unpublished work one has
tried to prove convergence under the Riemann hypothesis of subsequences of $\{F_n\}$ such as when $n$ is
restricted to run along the solutions of $\sum_{a=1}^n\mu(a)=0$. Another attempt by J. B. Conrey and G. Myerson
\cite{conrey} relates to a mollification of $F_n$, the\textit{ Selberg approximation}, defined in \cite{notes3} by

$$
S_n:=\sum_{a=1}^n\mu(a)\left(1-\frac{\log a}{\log n}\right)\rho_a.
$$
\ \\  
A common problem to these sequences is that if they converge at all to $-\chi$ in $\hilbert$ they must do so very
slowly: it is known \cite{notes3} that for any $F=\sum_{k=1}^n c_k \rho_{a_k}$, $a_k\geq1$, if $N=\max a_k$, then

\begin{equation}\label{slow}
\|F-\chi\|_\hilbert\geq \frac{C}{\sqrt{\log N}},
\end{equation}
\ \\
for an absolute constant C that has recently been sharpened by J. F. Burnol \cite{burnol}. This, as well as
considerations of
\textit{ summability} of series, led the author in
\cite{AB} as well as here to try to employ symultaneously, as it were, the whole range of $a\in[1,\infty)$. Thus we
define for complex
$s$ and
$x>0$ the functions

\begin{equation}\label{fs}
f_s(x):=\sum_{a=1}^\infty \frac{\mu(a)}{a^s}\rho_a(x).
\end{equation}
For fixed $x>0$ this is a meromorphic function of $s$ in the complex plane since

$$
f_s(x)=\frac{1}{x\zeta(s+1)}-\sum_{a\leq1/x}\frac{\mu(a)}{a^s}\left[\frac{1}{ax}\right],
$$ 
where the finite sum on the right is an entire function; thus $f_s$ is seen to be a sort of correction of
$1/\zeta(s)$. Assuming the Riemann hypothesis we shall prove for small positive $\epsilon$ that

$$
f_\epsilon\in\overline{\Bnat},
$$
and then,\textit{ unconditionally}, that

$$
f_\epsilon\hlimit -\chi,\ \ \ (\epsilon\downarrow0),
$$
so that $\chi\in\overline{\Bnat}$.

\section{The Proof}

\subsection{Two technical lemmas}
Here $s=\sigma+i\tau$ with $\sigma$ and $\tau$ real. The well-known theorem of Littlewood (see \cite{titchmarsh}
Theorem 14.25 (A)) to the effect that under the Riemann hypothesis $\sum_{a\geq1}\mu(a)a^{-s}$ converges to
$1/\zeta(s)$ for $\Re(s)>1/2$ has been provided in the more general setting of $\Re(s)>\alpha$ with a precise error
term by M. Balazard and E. Saias (\cite{notes1}, Lemme 2). We quote their lemma here for the sake of convenience.

\begin{lemma}[Balazard-Saias]\label{lemme2}
Let $1/2\leq\alpha<1$, $\delta>0$, and $\epsilon>0$. If $\zeta(s)$ does not vanish in the half-plane $\Re(s)>\alpha$,
then for $n\geq2$ and $\alpha+\delta\leq\Re(s)\leq1$ we have

\begin{equation}
\sum_{a=1}^n\frac{\mu(a)}{a^s}=
\frac{1}{\zeta(s)}+O_{\alpha,\delta,\epsilon}\left(n^{-\delta/3}(1+|\tau|)^\epsilon\right)
\end{equation} 
\end{lemma}
\ \\
It is important to note that the next lemma is independent of the Riemann or even the Lindel\"{o}f hypothesis.

\begin{lemma}\label{zratio}
For $0\leq\epsilon\leq \epsilon_0<1/2$ there is a positive constant $C=C(\epsilon_0)$ such that
for all
$\tau$
\begin{equation}\label{zratioeq}
\left|\frac{\zeta(\frac{1}{2}-\epsilon+i\tau)}{\zeta(\frac{1}{2}+\epsilon+i\tau)}\right|\leq
C\left(1+|\tau|\right)^\epsilon.
\end{equation}
\end{lemma}
\begin{proof}
We bring in the functional equation of $\zeta(s)$ to bear as follows
\begin{eqnarray}\nonumber
\left|\frac{\zeta(\frac{1}{2}-\epsilon+i\tau)}{\zeta(\frac{1}{2}+\epsilon+i\tau)}\right|&=&
\left|\frac{\zeta(\frac{1}{2}-\epsilon-i\tau)}{\zeta(\frac{1}{2}+\epsilon+i\tau)}\right|\\\nonumber
&=&\pi^{-\epsilon}\left|
\frac{\Gamma(\frac{1}{4}+\frac{1}{2}\epsilon+\frac{1}{2}i\tau)}{\Gamma(\frac{1}{4}-\frac{1}{2}\epsilon+\frac{1}{2}i\tau)}
\right|,
\end{eqnarray}
then the conclusion follows easily from well-known asymptotic formulae for the gamma function in a vertical strip
(\cite{rademacher} (21.51), (21.52)).
\end{proof}

\subsection{The proof proper of Theorem \ref{natconj}} 

It is clear that we need not prove the \textit{if} part of Theorem \ref{natconj}. So let us assume that the Riemann
hypothesis is true. We define

$$
f_{\epsilon,n}:=\sum_{a=1}^n\frac{\mu(a)}{a^\epsilon}\rho_a,\ \ \ (\epsilon>0).
$$ 
It is easy to see that

\begin{equation}\label{fepsilonn}
f_{\epsilon,n}(x)=\frac{1}{x}\sum_{a=1}^n\frac{\mu(a)}{a^{1+\epsilon}}
-\sum_{a=1}^n\frac{\mu(a)}{a^\epsilon}\left[\frac{1}{ax}\right],
\end{equation}
then, noting that the terms of the right-hand sum drop out when $a>1/x$, we obtain the pointwise limit

\begin{equation}\label{fepsilon1}
f_\epsilon(x)=\lim_{n\rightarrow\infty}f_{\epsilon,n}(x)=\frac{1}{x\zeta(1+\epsilon)}
-\sum_{a\leq 1/x}\frac{\mu(a)}{a^\epsilon}\left[\frac{1}{ax}\right].
\end{equation}
Then again for fixed $x>0$ we have
\begin{equation}\label{endlimit}
\lim_{\epsilon\downarrow0}f_\epsilon(x)=-\sum_{a\leq 1/x}\mu(a)\left[\frac{1}{ax}\right]=-\chi(x),
\end{equation}
by the fundamental property on M\"{o}bius numbers. The task at hand now is to prove these pointwise limits are also
valid in the $\hilbert$-norm. To this effect we introduce a new Hilbert space

$$
\kspace:=L_2((-\infty,\infty),(2\pi)^{-1/2}dt),
$$
\ \\
and note that by virtue of Plancherel's theorem the Fourier-Mellin map
$\M$ defined by

\begin{equation}\label{isometry}
\M(f)(\tau):=\int_0^\infty x^{-\frac{1}{2}+i\tau}f(x)dx,
\end{equation}
\ \\
is an invertible isometry from $\hilbert$ to $\kspace$. A well-known identity, which is at the root of the
Nyman-Beurling formulation, probably due to Titchmarsh (\cite{titchmarsh}, (2.1.5)), namely

$$
-\frac{\zeta(s)}{s}=\int_0^\infty x^{s-1}\rho_1(x)dx,\ \ \ (0<\Re(s)<1),
$$
immediately yields, denoting $X_{\epsilon}(x)=x^{-\epsilon}$,

\begin{equation}\label{mfepsilonn}
\M(X_\epsilon f_{2\epsilon,n})(\tau)=
-\frac{\zeta(\frac{1}{2}-\epsilon+i\tau)}{\frac{1}{2}-\epsilon+i\tau}
\sum_{a=1}^n \frac{\mu(a)}{a^{\frac{1}{2}+\epsilon+i\tau}},\ \ \ (0<\epsilon<1/2).
\end{equation}
\ \\
By a theorem of Littlewood (\cite{titchmarsh}, Theorem 14.25 (A)) if we let $n\rightarrow\infty$ in the right-hand side
of (\ref{mfepsilonn}) we get the pointwise limit

\begin{equation}\label{prelimit}
-\frac{\zeta(\frac{1}{2}-\epsilon+i\tau)}{\frac{1}{2}-\epsilon+i\tau}
\sum_{a=1}^n \frac{\mu(a)}{a^{\frac{1}{2}+\epsilon+i\tau}} \rightarrow
-\frac{\zeta(\frac{1}{2}-\epsilon+i\tau)}{\zeta(\frac{1}{2}+\epsilon+i\tau)}\frac{1}{\frac{1}{2}-\epsilon+i\tau}.
\end{equation}
\ \\
To see that this limit also takes place in $\hilbert$ we choose the parameters in Lemma \ref{lemme2} as
$\alpha=1/2$, $\delta=\epsilon>0$, $\epsilon\leq 1/2$, and $n\geq2$ to obtain

$$
\sum_{a=1}^n \frac{\mu(a)}{a^{\frac{1}{2}+\epsilon+i\tau}}=
\frac{1}{\zeta(\frac{1}{2}+\epsilon+i\tau)}+O_\epsilon\left((1+|\tau|)^\epsilon\right).
$$ 
If we now use Lemma \ref{zratio} and the Lindel\"{o}f hypothesis applied to the abscissa $1/2-\epsilon$, which follows
from the Riemann hypothesis, we obtain a positive constant $K_\epsilon$ such that for all real $\tau$

$$
\left|-\frac{\zeta(\frac{1}{2}-\epsilon+i\tau)}{\frac{1}{2}-\epsilon+i\tau}
\sum_{a=1}^n \frac{\mu(a)}{a^{\frac{1}{2}+\epsilon+i\tau}}\right|
\leq K_\epsilon(1+|\tau|)^{-1+2\epsilon}.
$$
\ \\
It is then clear that for $0<\epsilon\leq\epsilon_0<1/2$ the left-hand side of (\ref{prelimit}) is uniformly
majorized by a function in $\kspace$. Thus the convergence does take place in $\kspace$ which implies that

$$
X_\epsilon f_{2\epsilon,n}\hlimit X_\epsilon f_{2\epsilon}.
$$
\ \\
But $x^{-\epsilon}>1$ for $0<x<1$, and for $x>1$

\begin{equation}\label{forx>1}
f_{2\epsilon,n}(x)=\frac{1}{x}\sum_{a=1}^n\frac{\mu(a)}{a^{1+\epsilon}}\ll\frac{1}{x}, \ \ \ (x>1),
\end{equation}
uniform in $n$, which easily implies that one also has $\hilbert$-convergence for $f_{2\epsilon,n}$ as
$n\rightarrow\infty$. The factor $2$ in the subindex is unessential, so that we now have for sufficiently small
$\epsilon>0$ that

$$
f_{\epsilon,n}\hlimit f_\epsilon\in\overline{\Bnat},
$$
\ \\  
as was announced above. Moreover, since we have identified the pointwise limit in (\ref{prelimit}) we now
have

$$
\M(X_\epsilon f_{2\epsilon})(t)=
-\frac{\zeta(\frac{1}{2}-\epsilon+i\tau)}{\zeta(\frac{1}{2}+\epsilon+i\tau)}\frac{1}{\frac{1}{2}-\epsilon+i\tau}.
$$
\ \\
Now we apply Lemma \ref{zratio} and obtain, at this juncture now \textit{without the assumption of the Riemann
hypothesis}, that
$\M(X_\epsilon f_{2\epsilon})$ converges in $\kspace$, thus $X_\epsilon f_{2\epsilon}$ converges in $\hilbert$, and
this means that $f_\epsilon$ also converges in $\hilbert$ as $\epsilon\downarrow0$ by an argument entirely similar to
that used for $f_{\epsilon,n}$. The identification of the pointwise limit in (\ref{endlimit}) finally gives

$$
f_\epsilon\hlimit -\chi,
$$   
\ \\
which concludes the proof.

\section{A corollary and a constructive approximating sequence}
The proof of Theorem \ref{natconj} provides in turn a new proof, albeit of a stronger theorem, of the
Nyman-Beurling criterion which bypasses the deep and complicated Hardy space techniques. One should extend it to the
$L_p$ case, that is, to the condition that $\zeta(s)$ does not vanish in a half-plane $\Re(s)>1/p$. It should be clear
also that we have shown this special equivalence criterion to be true:
\begin{corollary}
The Riemann hypothesis is equivalent to the $\hilbert$-convergence of $f_{\epsilon,n}$ as
$n\rightarrow\infty$ for all sufficiently small $\epsilon>0$. 
\end{corollary}
In essence what has been done is to apply a summability method to the old natural approximation. The convergence on
special subsequences both of $n$ and of $\epsilon$ is also necessary and sufficient, and it is proposed here to look
further into this alongside with other summability methods for the natural approximation. Indeed, it is worth noticing
that E. Bombieri\footnote{Personal communication 20 Feb. 2002.} indicated to us some advantages of using a Ces\`{a}ro
smoothing of the $f_{\epsilon,n}$, that is, to substitute instead the
functions
\begin{equation}\label{eb}
\sum_{a=1}^n\frac{\mu(a)}{a^\epsilon}\left(1-\frac{a}{n}\right)\rho_a,
\end{equation}
which clearly deserves closer study. Independently of each other M. Balazard and E. Bombieri  mentioned
to the author that a suitable choice of $\epsilon$ would lead to a
quantitative estimate of the error of approximation as some inverse
power of $\log\log n$. Note that we did not employ the dependence on $n$ in the Balazard-Saias Lemma
\ref{lemme2}. This dependence is obviously connected to the slowness of approximation to $\chi$ indicated
in (\ref{slow}). We owe special thanks to M. Balazard for the statement and proof of the following proposition:

\begin{proposition}[M. Balazard]\label{michel}
If the Riemann hypothesis is true then there is a constant $c>0$ such that the distance in $\hilbert$ between $\chi$ and
\begin{equation}\label{approx}
-\sum_{a=1}^n \mu(a)e^{-c\frac{\log a}{\log\log n}}\rho_a
\end{equation}
is $\ll (\log\log n)^{-1/3}$. 
\end{proposition}

\begin{proof}
We shall only sketch the proof of this proposition. Applying the Fourier-Mellin map (\ref{isometry}) to
$f_{\epsilon,n}+\chi$ we have from Plancherel's theorem that 

\begin{eqnarray}\nonumber
2\pi\|f_{\epsilon,n}+\chi\|_\hilbert^2&=&\int_{\Re(z)=1/2}
\left|\zeta(z)\sum_{a=1}^n\frac{\mu(a)}{a^{z+\epsilon}}-1\right|^2 \frac{|dz|}{|z|^2}\\\nonumber
 &\leq&2\int_{\Re(z)=1/2}
\left|\zeta(z)\left(\sum_{a=1}^n\frac{\mu(a)}{a^{z+\epsilon}}-\frac{1}{\zeta(z+\epsilon)}\right)\right|^2
\frac{|dz|}{|z|^2}\\\nonumber
& &+
2\int_{\Re(z)=1/2}\left|\frac{\zeta(z)}{\zeta(z+\epsilon)}-1\right|^2\frac{|dz|}{|z|^2}.
\end{eqnarray}
The second integral on the right-hand side above is estimated to be $\ll\epsilon^{2/3}$ as follows. If the distance
between $z= \frac{1}{2}+it$ and the nearest zero of $\zeta$ is larger than $\delta$, say, the upper bound

$$
\left \lvert \frac{\zeta (z)}{\zeta(z+\epsilon)} -1 \right \rvert \ll \epsilon \delta^{-1}  (|t| +1)^{1/4}
$$
follows from the classical estimate

$$
\frac{\zeta' (s)}{\zeta(s)} = \sum_{|\gamma - \tau | \leq 1} \frac{1}{s-\rho} + O \bigl ( \log (2 + |\tau|) \bigr ),
$$
where $s=\sigma + i \tau$, $1/2 \leq \sigma \leq 3/4$, $\tau \in \Real$ and $\rho =\beta + i \gamma$ denotes a generic
zero of the $\zeta$ function, by integration and exponentiation, provided $\epsilon/\delta$ is small enough. In the other
case, one uses an estimate of Burnol \cite{burnolineq} stating that under the Riemann hypothesis

\begin{equation}\label{burnolestimate}
\left|\frac{\zeta(z)}{\zeta(z+\epsilon)}\right|\ll |z|^{\epsilon/2}, \ \ \Re(z)=1/2, \ \ 0<\epsilon\leq 1/2.
\end{equation}
Integrating these two inequalities on the corresponding sets, one gets an upper bound $\ll \epsilon^2/\delta^2 +
\delta$, and chooses $\delta = \epsilon^{2/3}$. 
The first integral, on the other hand, succumbs to a special form of the Balazard-Saias Lemma \ref{lemme2}, namely

$$
\sum_{a=1}^n\frac{\mu(a)}{a^{\frac{1}{2}+\epsilon+it}} = \frac{1}{\zeta(\frac{1}{2}+\epsilon+it)}+
O(n^{-\epsilon/3}e^{b\mathcal L(t)}), \ \ \ c/\log\log n \leq \epsilon\leq1/2, 
$$
where $\mathcal{L}(t):=\log(|t|+3)/\log\log(|t|+3)$. We thus have

\begin{eqnarray}\nonumber
2\pi\|f_{\epsilon,n}+\chi\|_\hilbert^2 &\ll& n^{-2\epsilon/3}\int_{-\infty}^\infty e^{O(\mathcal{L}(t))}
\frac{dt}{\frac{1}{4}+t^2}+\epsilon^{2/3}\\\nonumber
 &\ll& n^{-2\epsilon/3}+\epsilon^{2/3},
\end{eqnarray}
provided that $\epsilon\geq c/\log\log n$, whereupon one chooses $\epsilon=c/\log\log n$ to reach the conclusion.
\end{proof}
\ \\
\ \\
\textbf{Acknowledgements.} We are happy to thank E. Bombieri for thorough comments and suggestions
(in particular (\ref{eb})), to M. Balazard and E. Saias for pointing out their lemma \ref{lemme2}, to M. Balazard for
Proposition \ref{michel} estimating the error term for the approximation (\ref{approx}), to J. F. Burnol [8], who, after
the initial version of this paper, provided an interesting approach centered on the zeta estimate (\ref{burnolestimate})
which turned out also to be used in the Balazard Proposition \ref{michel}.

\bibliographystyle{amsplain}

\ \\
\noindent Luis B\'{a}ez-Duarte\\
Departamento de Matem\'{a}ticas\\
Instituto Venezolano de Investigaciones Cient\'{\i}ficas\\
Apartado 21827, Caracas 1020-A\\
Venezuela\\

\end{document}